\documentclass[11pt,twoside]{amsart}
\usepackage{amsmath, amsthm, amscd, amsfonts, amssymb, graphicx, color}
\usepackage[bookmarksnumbered, colorlinks, plainpages]{hyperref}
\usepackage{pgf,tikz}
\usetikzlibrary{arrows}
\usepackage[utf8]{inputenc}
\usepackage{pstricks-add}
\input{mathrsfs.sty}
\newtheorem{thm}{Theorem}[section]
\newtheorem{cor}[thm]{Corollary}

\newtheorem{defn}[thm]{Definition}
\newtheorem{rem}[thm]{\bf{Remark}}
\newtheorem{example}[thm]{\bf{Example}}

\numberwithin{equation}{section}

\begin{document}

%------------------------------------------------------------------------------------%

\title{A generalization of Hall's theorem for $k$-uniform $k$-partite hypergraphs}
\author{Reza jafarpour-Golzari}
\address{Department of Mathematics, Institute for Advanced Studies
in Basic Science (IASBS), P.O.Box 45195-1159, Zanjan, Iran}

\email{r.golzary@iasbs.ac.ir}

%\address{Department of Mathematics, Institute for Advanced Studies
%in Basic Science (IASBS), P.O.Box 45195-1159, Zanjan, Iran}

%\email{rashidzn@iasbs.ac.ir}

\thanks{{\scriptsize
\hskip -0.4 true cm MSC(2010): Primary: 05E40; Secondary: 05C65, 05D15.
\\
%$*$Corresponding author
\newline\indent{\scriptsize}}}

\maketitle

%------------------------------------------------------------------------------------%

\begin{abstract}
In this paper we prove a generalized version of Hall's theorem for hypergraphs. More precisely, let $\mathcal{H}$ be a $k$-uniform $k$-partite hypergraph with some ordering on parts as $V_{1}, V_{2}, \ldots , V_{k}$. such that the subhypergraph generated on $\bigcup_{i=1}^{k-1}V_{i}$ has a unique perfect  matching. In this case, we give a necessary and sufficient condition for having a matching of size $t=|V_{1}|$ in $\mathcal{H}$. Some relevant results and counterexamples are given as well. 
\end{abstract}

\vskip 0.2 true cm

%------------------------------------------------------------------------------------%

\section{\bf Introduction}
\vskip 0.4 true cm

We refer to \cite{Bon} and \cite{Ber} for elementary backgrounds in graph and hypergraph theory respectively.

Let $G$ be a simple finite graph with vertex set $V(G)$ and edge set $E(G)$. A matching in $G$, is a set $M$ of pairwise disjoint edges of $G$. A matching $M$ is said to be a perfect matching, if every $x\in V(G)$ lies in one of elements of $M$. A matching $M$ in $G$, is maximum whenever for every matching $M^{'}$, $|M^{'}|\leq |M|$.

For every set of vertices $A$, $N(A)$ which is called the neighborhood of $A$ is the set of vertices which are adjacent with at least one element of $A$. The following theorem is known as Hall's theorem in bipartite graphs.

\begin{thm}
(\cite{Bon} Theorem 5.2) Let $G$ be a bipartite graph with bipartition $(X, Y)$. Then $G$ contains a matching that saturates every vertex in $X$, if and only if
\[|N(S)|\geq |S| \ \ \ \ for \ all \ S\subseteq X.\]
\end{thm}

 A vertex cover in $G$, is a subset $C$ of $V(G)$ such that for every edge $e$ of $G$, $e$ intersects $C$. A vertex cover $C$ is called a minimum vertex cover, if for every vertex cover $C^{'}$, $|C|\leq |C^{'}|$. The following theorem is known as K\"{o}nig's theorem in graph theory.

\begin{thm}
(\cite{Bon} Theorem 5.3) In a bipartite graph, the number of edges in a maximum matching is equal to the number of vertices in a minimum vertex cover.
\end{thm}

Let $V$ be a finite nonempty set. A hypergraph $\mathcal{H}$ on $V$ is a collection of nonempty subsets of $V$ such that $\bigcup_{e\in\mathcal{H}}e=V$. Each subset is said to be a hyperedge and each element of $V$ is called a vertex. We denote the set of vertices and hyperedges of $\mathcal{H}$ by $V(\mathcal{H})$ and $E(\mathcal{H})$, respectively. Two vertices $x, y$ of a hypergraph are said to be adjacent whenever they lie in a hyperedge. 

A matching in the hypergraph $\mathcal{H}$ is a set $M$ of pairwise disjoint hyperedges of $\mathcal{H}$. A perfect matching is a matching such that every $x\in V(\mathcal{H})$ lies in one of its elements. A matching $M$ in $\mathcal{H}$ is called a maximum matching whenever for every matching $M^{'}$, $|M^{'}|\leq |M|$. 

In a hypergraph $\mathcal{H}$, a subset $C$ of $V(\mathcal{H})$ is called a vertex cover if every hyperedge of $\mathcal{H}$ intersects $C$. A vertex cover $C$ is said to be minimum if for every vertex cover $C^{'}$, $|C|\leq |C^{'}|$. We denote the number of hyperedges in a maximum matching of the hypergraph $\mathcal{H}$  by $\alpha^{'}(\mathcal{H})$ and the number of vertices in a minimum vertex cover of $\mathcal{H}$ by $\beta(\mathcal{H})$. 

A hypergraph $\mathcal{H}$ is said to be simple or a clutter if non of its two distinct hyperedges contains another. A hypergraph is called $t$-uniform (or $t$-graph), if all its hyperedges have the same size $t$. A hypergraph $\mathcal{H}$ is said to be $r$-partite ($r\geq 2$), whenever $V(\mathcal{H})$ can be partitioned to $r$ subsets such that for every two vertices $x, y$ in one part, $x$ and $y$ are not adjacent. If $r=2, 3$, the hypergraph is said to be bipartite and tripartite respectively. 

Several researches have been done about matching and existence of perfect matching in hypergraphs (see for instance \cite{Aha 1}, \cite{F}, \cite{K}). Also some attemps have been produced in generalization of Hall's theorem and K\"{o}nig's theorem to hypergraphs (see \cite{Aha 2}, \cite{Aha 3}, \cite{Aha 4}, \cite{Aha 5}, \cite{Hax 1}, \cite{Hax 2}). 

\begin{defn}
Let $\mathcal{H}$ be a $k$-uniform hypergraph with $k\geq 2$. A subset $\mathfrak{e}\subseteq V(\mathcal{H})$ of size $k-1$ is called a submaximal edge if there is a hyperedge containing $\mathfrak{e}$. For a submaximal edge $\mathfrak{e}$, define the neighborhood of $\mathfrak{e}$ as the set $N(\mathfrak{e}):=\{v\in V(\mathcal{H})| \ \mathfrak{e}\cup \{v\}\in E(\mathcal{H})\}$. 
\end{defn}

For a set $A$ consisting of submaximal edges of $\mathcal{H}$, $\{v\in V(\mathcal{H})| \ \exists \mathfrak{e}\in A, v\in N(\mathfrak{e})\}$ is denoted by $N(A)$.

\begin{defn}
Let $\mathcal{H}$ be a hypergraph, and $\emptyset\neq V^{'}\subseteq V(\mathcal{H})$. The subhypergraph generated on $V^{'}$ is
\[<V^{'}>:= \{e\cap V^{'}| \ e\in E(\mathcal{H}), e\cap V^{'}\neq\emptyset\}.\]  
\end{defn}

If $\mathcal{H}$ is a $k$-uniform $k$-partite hypergraph with parts $V_{1}, V_{2}, \ldots , V_{k}$, it is clear that the subhypergraph generated on the union of every $k-1$ distinct parts is a $(k-1)$-uniform $(k-1)$-partite hypergraph.

Let $\mathfrak{A}=(A_{1}, \ldots , A_{n})$ be a family of subsets of a set $E$. A subset $\{x_{1}, \ldots , x_{n}\}_{\neq}$ of $E$ is said to be a transversal (or SDR) for $\mathfrak{A}$, if for every $i \ (1\leq i\leq n), \ x_{i}\in A_{i}$. A partial transversal (partial SDR) of length $l \ (1\leq l\leq n-1)$ for $\mathfrak{A}$, is a transversal for a subfamily of $\mathfrak{A}$ with $l$ sets.\cite{Bry}

The following theorem is known as Hall's theorem in combinatorics.

\begin{thm}
(\cite{Bry} Theorem 4.1) The family $\mathfrak{A}=(A_{1}, \ldots , A_{n})$ of
 subsets of a set $E$ has a transversal if and only if
 \[|\bigcup_{i\in I^{'}} A_{i}|\geq |I^{'}|, \ \ \ \forall I^{'}\subseteq \{1, \ldots , n\}.\]
\end{thm}

\begin{cor}
(\cite{Bry} Corollary 4.3) The family $\mathfrak{A}=(A_{1}, \ldots , A_{n})$ of
 subsets of a set $E$ has a partial transversal of length $l(>0)$ if and only if
 \[|\bigcup_{i\in I^{'}} A_{i}|\geq |I^{'}|-n+l, \ \ \ \forall I^{'}\subseteq \{1, \ldots , n\}.\]
\end{cor}
%------------------------------------------------------------------------------------%

\section{\bf {\bf \em{\bf The main results}}}
\vskip 0.4 true cm

Now we are ready to present our first theorem.

\begin{thm}
Let $\mathcal{H}$ is a $k$-uniform $k$-partite hypergraph with some ordering on parts, as $V_{1}, V_{2}, \ldots , V_{k}$ such that the subhypergraph generated on $\bigcup_{i=1}^{k-1}V_{i}$ has a unique perfect matching $M$. Then $\mathcal{H}$ has a matching of size $t=|V_{1}|$, if and only if for every subset $A$ of $M$, $|N(A)|\geq |A|$. 
\begin{proof}
Let $t=|V_{1}|$ and let the elements of $M$ are $\mathfrak{e}_{1}, \ldots , \mathfrak{e}_{t}$. Let $\mathcal{H}$ has a matching of size $t$ with elements $e_{1}, \ldots , e_{t}$. By uniqueness of $M$, $M=\{e_{1}-V_{k}, \ldots , e_{t}-V_{k}\}$. Therefore
\[(N(\mathfrak{e}_{1}), \ldots , N(\mathfrak{e}_{t}))=(N(e_{i_{1}}-V_{k}), \ldots , N(e_{i_{t}}-V_{k})).\]
Then the family $(N(\mathfrak{e}_{1}), \ldots , N(\mathfrak{e}_{t}))$ has an SDR. Then by Theorem 1.3
\[|\bigcup_{i\in I} N(\mathfrak{e}_{i})|\geq |I|, \ \ \ \forall I\subseteq \{1, \ldots , t\}\]
and therefore for every subset $A$ of $M$, $|N(A)|\geq |A|$.

Conversely let for every subset $A$ of $M$, we have $|N(A)|\geq |A|$. Now $(N(\mathfrak{e}_{1}), \ldots , N(\mathfrak{e}_{t}))$ is a family such that 
\[|\bigcup_{i\in I} N(\mathfrak{e}_{i})|\geq |I|, \ \ \ \forall I\subseteq \{1, \ldots , t\}.\]
Therefore by Theorem 1.3, the mentioned family has an SDR. That is, there are distinct elements $x_{1}, \ldots , x_{t}$ of $V_{k}$ such that $x_{j}\in N(\mathfrak{e}_{j})$. Now for every $1\leq j\leq t$, $\mathfrak{e}_{j}\cup\{x_{j}\}$ is a hyperedge of $\mathcal{H}$ and these hyperedges are pairwise disjoint. Then they form a matching of size $t$ for $\mathcal{H}$.
\end{proof}
\end{thm}

\begin{cor}
Let $\mathcal{H}$ be a $k$-uniform $k$-partite hypergraph with some ordering on parts as $V_{1}, V_{2}, \ldots , V_{k}$ where $|V_{1}|=|V_{2}|= \cdots =|V_{k}|$, such that the subhypergraph generated on $\bigcup_{i=1}^{k-1}V_{i}$ has a unique perfect matching $M$. Then $\mathcal{H}$ has a perfect matching if and only if for every subset $A$ of $M$, $|N(A)|\geqslant |A|$.  
\end{cor}

\begin{rem}
Theorem 2.1 implies Theorem 1.1 (Hall's theorem) in case $k=2$.
\end{rem}

\begin{rem}
In Theorem 2.1, if the hypothesis of uniqueness of perfect matching of subhypergraph generated on $\bigcup_{i=1}^{k-1}V_{i}$ is removed, only one side of theorem will remains correct. That is, from this fact that for every subset $A$ of $M$, $|N(A)|\geqslant |A|$, we conclude that $\mathcal{H}$ has a matching of size $t=|V_{1}|$. The following example shows that the inverse case is not true in general.
\end{rem}

\begin{example}\rm
Assume the 3-uniform 3-partite hypergraph $\mathcal{H}$ with the following presentation.

%\begin{figure}
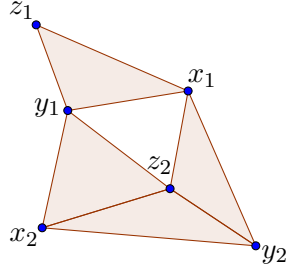
\begin{center}
\definecolor{zzttqq}{rgb}{0.6,0.2,0.}
\definecolor{qqqqff}{rgb}{0.,0.,1.}
\begin{tikzpicture}[line cap=round,line join=round,>=triangle 45,x=1.0cm,y=1.0cm]
\clip(8.,2.) rectangle (13.4,6.9);
\fill[color=zzttqq,fill=zzttqq,fill opacity=0.1] (11.26,5.1) -- (9.24,5.98) -- (9.66,4.84) -- cycle;
\fill[color=zzttqq,fill=zzttqq,fill opacity=0.1] (11.26,5.1) -- (12.16,3.04) -- (11.02,3.8) -- cycle;
\fill[color=zzttqq,fill=zzttqq,fill opacity=0.1] (11.02,3.8) -- (9.66,4.84) -- (9.32,3.28) -- cycle;
\fill[color=zzttqq,fill=zzttqq,fill opacity=0.1] (12.16,3.04) -- (9.32,3.28) -- (11.02,3.8) -- cycle;
\draw [color=zzttqq] (11.26,5.1)-- (9.24,5.98);
\draw [color=zzttqq] (9.24,5.98)-- (9.66,4.84);
\draw [color=zzttqq] (9.66,4.84)-- (11.26,5.1);
\draw [color=zzttqq] (11.26,5.1)-- (12.16,3.04);
\draw [color=zzttqq] (12.16,3.04)-- (11.02,3.8);
\draw [color=zzttqq] (11.02,3.8)-- (11.26,5.1);
\draw [color=zzttqq] (11.02,3.8)-- (9.66,4.84);
\draw [color=zzttqq] (9.66,4.84)-- (9.32,3.28);
\draw [color=zzttqq] (9.32,3.28)-- (11.02,3.8);
\draw [color=zzttqq] (12.16,3.04)-- (9.32,3.28);
\draw [color=zzttqq] (9.32,3.28)-- (11.02,3.8);
\draw [color=zzttqq] (11.02,3.8)-- (12.16,3.04);
\draw (8.74,6.42) node[anchor=north west] {$z_{1}$};
\draw (9.08,5.14) node[anchor=north west] {$y_{1}$};
\draw (8.75,3.4) node[anchor=north west] {$x_{2}$};
\draw (10.57,4.36) node[anchor=north west] {$z_{2}$};
\draw (12.1,3.18) node[anchor=north west] {$y_{2}$};
\draw (11.12,5.52) node[anchor=north west] {$x_{1}$};
\begin{scriptsize}
\draw [fill=qqqqff] (11.26,5.1) circle (1.5pt);
\draw [fill=qqqqff] (9.24,5.98) circle (1.5pt);
\draw [fill=qqqqff] (9.66,4.84) circle (1.5pt);
\draw [fill=qqqqff] (12.16,3.04) circle (1.5pt);
\draw [fill=qqqqff] (11.02,3.8) circle (1.5pt);
\draw [fill=qqqqff] (9.32,3.28) circle (1.5pt);
\end{scriptsize}
\end{tikzpicture}
% \caption{}
  \end{center}
  \label{graph 1}
% \end{figure}
Indeed, $\mathcal{H}=\{\{x_{1}, y_{1}, z_{1}\}, \{x_{1}, y_{2}, z_{2}\}, \{x_{2}, y_{2}, z_{2}\}, \{x_{2}, y_{1}, z_{2}\}\}$ where the parts of $\mathcal{H}$ are
\[V_{1}=\{x_{1}, x_{2}\}, \ V_{2}=\{y_{1}, y_{2}\}, \ V_{3}=\{z_{1}, z_{2}\}.\] 
In this case there is a perfect matching $M_{1}=\{\{x_{2}, y_{1}\}, \{x_{1}, y_{2}\}\}$
for subhypergraph generated on $V_{1}\cup V_{2}$. Although the hypergraph $\mathcal{H}$ has a matching $M^{'}=\{\{x_{1}, y_{1}, z_{1}\}, \{x_{2}, y_{2}, z_{2}\}\}$ of size 2, if $A=M_{1}$, we have $N(A)=\{z_{2}\}$. Therefore $|N(A)|\ngeqslant |A|$. Note that $M_{1}$ is not the unique perfect matching of subhypergraph generated on $V_{1}\cup V_{2}$ because $M_{2}=\{\{x_{1}, y_{1}\}, \{x_{2}, y_{2}\}\}$ is also yet.
\end{example}

\begin{thm}
Let $\mathcal{H}$ be a $k$-uniform $k$-partite hypergraph with some ordering on parts as $V_{1}, V_{2}, \ldots , V_{k}$ such that the subhypergraph generated on $\bigcup_{i=1}^{k-1}V_{i}$ has a perfect matching $M$. If for every subset $A$ of $M$, we have $|N(A)|\geqslant |A|-p$ where $p$ is a fix integer and $1\leq p\leq t-1$, then $\mathcal{H}$ has a matching of size $t-p$, where $t$ is the size of $V_{1}$.
\begin{proof}
Let the elements of $M$ be $\mathfrak{e}_{1}, \ldots , \mathfrak{e}_{t}$. $(N(\mathfrak{e}_{1}), \ldots , N(\mathfrak{e}_{t}))$ is a family such that the cardinality of the union of each $s$ terms is greater than or equal to $s-t+(t-p)$. Then by Corollary 1.4, the family $(N(\mathfrak{e}_{1}), \ldots , N(\mathfrak{e}_{t}))$ has a partial SDR of size $t-p$. That is, there are distinct elements $y_{1}, \ldots , y_{t-p}$ of $V_{k}$ such that $y_{j}\in N(\mathfrak{e}_{i_{j}})$. Now for every $1\leq j\leq t-p$, $\mathfrak{e}_{i_{j}}\cup\{y_{j}\}$ is a hyperedge of $\mathcal{H}$ and these hyperedges are pairwise disjoint. Then they form a matching of size $t-p$ for $\mathcal{H}$. 
\end{proof}
\end{thm}

\begin{thm}
Let $\mathcal{H}$ be a $k$-uniform $k$-partite hypergraph with some ordering on parts as $V_{1}, V_{2}, \ldots , V_{k}$, and let $t=|V_{1}|$. Then $\mathcal{H}$ has a matching of size $t$ if and only if $\alpha^{'}=\beta=t$.
\begin{proof}
Let $\mathcal{H}$ has a matching of size $t$. We show that  $\alpha^{'}=\beta=t$. Clearly $\beta\geq\alpha^{'}$ because for covering each hyperedge of maximum matching, one vertex is needed. But since there is a matching of size $t$, then $\alpha^{'}\geq t$. Now $V_{1}$ is a minimal vertex cover of $\mathcal{H}$ because each hyperedge has only one vertex in $V_{1}$ and each vertex of $V_{1}$ lies in a hyperedge. Therefore $t\geq\beta$ which implies that $\alpha^{'}\geq\beta$. Then $\alpha^{'}=\beta$. The matching of size $t$ is the maximum matching because it covers all vertices of $V_{1}$.

Conversely, if $\alpha^{'}=\beta=t$, it is clear that $\mathcal{H}$ has a matching of size $t$.  
\end{proof}
\end{thm}

The following example shows that removing the condition $t=|V_{1}|$ in Theorem 2.7 is not possible even if the subhypergraph generated on union of every $k-1$ parts, has a perfect matching.

\begin{example}\rm
Assume 3-uniform 3-partite hypergraph $\mathcal{H}$ with the following presentation, where the parts of $\mathcal{H}$ are
\[V_{1}=\{1, 2\}, \ V_{2}=\{3, 4\}, \ V_{3}=\{5, 6\}.\]
Indeed $\mathcal{H}=\{\{1, 3, 5\}, \{2, 3, 6\}, \{2,4,5\}\}$.

%\begin{figure}
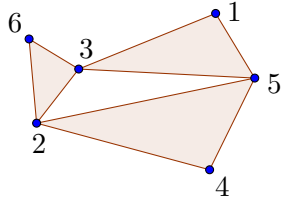
\begin{center}
\definecolor{zzttqq}{rgb}{0.6,0.2,0.}
\definecolor{qqqqff}{rgb}{0.,0.,1.}
\begin{tikzpicture}[line cap=round,line join=round,>=triangle 45,x=1.0cm,y=1.0cm]
\clip(6.,2.) rectangle (11.18,6.42);
\fill[color=zzttqq,fill=zzttqq,fill opacity=0.1] (7.78,4.72) -- (7.12,5.12) -- (7.22,4.) -- cycle;
\fill[color=zzttqq,fill=zzttqq,fill opacity=0.1] (7.78,4.72) -- (9.6,5.46) -- (10.12,4.6) -- cycle;
\fill[color=zzttqq,fill=zzttqq,fill opacity=0.1] (7.22,4.) -- (10.12,4.6) -- (9.52,3.38) -- cycle;
\draw (6.7,5.6) node[anchor=north west] {6};
\draw (7.02,3.98) node[anchor=north west] {2};
\draw (9.46,3.42) node[anchor=north west] {4};
\draw (7.65,5.25) node[anchor=north west] {3};
\draw (9.62,5.74) node[anchor=north west] {1};
\draw (10.15,4.78) node[anchor=north west] {5};
\draw [color=zzttqq] (7.78,4.72)-- (7.12,5.12);
\draw [color=zzttqq] (7.12,5.12)-- (7.22,4.);
\draw [color=zzttqq] (7.22,4.)-- (7.78,4.72);
\draw [color=zzttqq] (7.78,4.72)-- (9.6,5.46);
\draw [color=zzttqq] (9.6,5.46)-- (10.12,4.6);
\draw [color=zzttqq] (10.12,4.6)-- (7.78,4.72);
\draw [color=zzttqq] (7.22,4.)-- (10.12,4.6);
\draw [color=zzttqq] (10.12,4.6)-- (9.52,3.38);
\draw [color=zzttqq] (9.52,3.38)-- (7.22,4.);
\begin{scriptsize}
\draw [fill=qqqqff] (7.78,4.72) circle (1.5pt);
\draw [fill=qqqqff] (7.12,5.12) circle (1.5pt);
\draw [fill=qqqqff] (7.22,4.) circle (1.5pt);
\draw [fill=qqqqff] (9.6,5.46) circle (1.5pt);
\draw [fill=qqqqff] (10.12,4.6) circle (1.5pt);
\draw [fill=qqqqff] (9.52,3.38) circle (1.5pt);
\end{scriptsize}
\end{tikzpicture}
% \caption{}
  \end{center}
  \label{graph 2}
% \end{figure}
In this hypergraph we have the matching $\{\{1, 3, 5\}\}$ of size 1. But $\alpha^{'}\neq\beta$ because $\alpha^{'}=1$ and $\beta=2$. Note that each one of subhypergraph generated on $V_{1}\cup V_{2}$, $V_{2}\cup V_{3}$ and $V_{1}\cup V_{3}$ have a perfect matching.
\end{example}

%------------------------------------------------------------------------------------%

\vskip 0.4 true cm

%\begin{center}{\textbf{Acknowledgments}}
%\end{center}

%\vskip 0.4 true cm

%------------------------------------------------------------------------------------%

\end{document}